\documentclass{article}
\usepackage{mathrsfs}
\usepackage{CJK}
\usepackage{graphics}
\usepackage{graphicx}
\usepackage{epsfig}
\usepackage{amsmath}
\usepackage{amsthm}
\usepackage{amsfonts}
\usepackage{amssymb}
\usepackage{float}
\usepackage{CJK}
\usepackage{graphicx}
\usepackage{epstopdf}
\usepackage{amsthm}
\usepackage{color}
\usepackage{setspace}
\usepackage{geometry}
\usepackage{lineno}
\let\oldequation\eqaution
\let\oldequation\endeqaution

\begin{document}

\pagewiselinenumbers

\begin{spacing}{1.3}
\begin{CJK*}{GBK}{song}
\renewcommand\figurename{Figure}
\CJKtilde
\title{\bf The average order of a connected vertex set in $K_m \times P_n$}
\author{Mingyuan Ma\footnote{E-mail: 623mmy@sina.com.},\quad  Han Ren\footnote{E-mail: hren@math.ecnu.edu.cn}\\
\small{Army Engineering University of PLA, Nanjing,210001,China}\\
\small{School of Mathematical Sciences,East China Normal University,Shanghai,200241,China}
}
\date{}
\maketitle

{\bf Abstract.} Let $G$ be a connected graph. Let $N(G)$ and $S(G)$ be the number of connected sets of $G$ and the 
sum of the orders of these connected sets of $G$, respectively. Then $A(G)=\frac{S(G)}{N(G)}$ is called the average order of a connected set of $G$.
In this paper, we derive a closed-form formula for $A(K_m \times P_n)$, where $K_m \times P_n$ is the 
Cartesian product of the complete graph $K_m$ and the path $P_n$. 
\par

{\bf Keywords.} connected sets, average order, density of connected sets \par
\vspace{4mm}

\section{Introduction.}

All graphs in this paper are simple and finite. Let $G$ be a connected graph of $n$ vertices. A vertex subset $M \subseteq V(G)$ is called a \emph{connected set} of $G$ if the subgraph induced by $M$ is connected. Let $N(G)$ denote the number of connected sets in $G$. Let $\{M_1, M_2, \dots, M_{N(G)}\}$ be the collection of all connected sets of $G$. Then we denote $S(G) = \sum\limits^{N(G)}_{i=1}|M_i| $ as the sum of orders of the connected sets in $G$. Furthermore, let $$A(G) = \frac{S(G)}{N(G)} \quad  \text{and} \quad D(G) = \frac{A(G)}{n}$$ denote the average order of connected sets of $G$ and the \emph{density} of connected sets of vertices, respectively.

There are numerous research concerning the average order and density of connected sets of a graph. In 1983, Jamison[6] proved that among all trees of order n, the path $P_n$ minimizes the average order of a subtree. In 2010, Vince and Wang[16] proved that if $T$ is a tree all of whose internal vertices have degree at least three, then $\frac{1}{2}\leq D(T) < \frac{3}{4}$. In 2018, Kroeker, Mol and Oellermann[8] proved that $\frac{n}{2}\leq A(G) \leq\frac{n+1}{2}$ for a connected graph of order $n$.
For ladders and circular ladders, Vince[14] gave explicit closed formulas for both the number and the average order of connected induced
subgraphs of them in terms of the classic Pell numbers in 2021. For example, the average order of a connected set of the ladder $L_n$ is $$A(L_n)=\frac{(32-45\bar\beta(n)-32\beta(n))+n(10+21\beta(n)+30\bar\beta(n))}{2(\beta(n+3)-4n-7)},\eqno(1.1)$$
 where $\beta(n)$ denotes Pell-Lucas number and $\bar{\beta}(n)$ denotes Pell number.

Let $H$ and $Q$ be two graphs with vertex sets $V(H)=\{u_1, u_2, \ldots, u_s\}$ and $V(Q)=\{v_1, v_2, \ldots, v_t\}$, respectively.
The \emph{Cartesian product} of $H$ and $Q$, denoted by $H\times Q$, is the graph with vertex set 
$V(H)\times V(Q)$ and edge set $\{(u_i,v_p)(u_j,v_q)|u_iu_j\in E(H)\, \mbox{and}\, v_p=v_q,\,\mbox{or}\, v_pv_q\in E(Q)\,  \mbox{and}\,  u_i=u_j \}$.
Let $P_l$ be the path of $l$ vertices. In fact, $L_n$ is exactly $P_n\times P_2$. 
In 2021, Vince [14] proposed the following question: 
Find a formula for $N(P_n\times P_n)$. 
So far as we know,
there is not any formula for $N(P_n\times P_n)$ for an arbitrary $n$. In 2025, Ma et al.[10] 
obtained a formula for $N(K_m\times P_n)$, where $K_m$ is the complete graph of $m$ vertices.
However, they do not give a formula for $A(K_m\times P_n)$. In this paper we will establish 
a closed-form formula for $A(K_m\times P_n)$.

The arrangement of the paper is as follows:
In Section 2, we introduce some notations. In particular, we introduce the recurrence matrix $A_m$.
In Section 3, we develop a matrix method for calculating the average order of connected subsets in $K_m\times P_n$ and provide a general formula together with its proof. Subsequently, we give a formula for the density of connected sets in $K_m\times P_n$.
In Section 4, we use the formulas from Section 3 to calculate $A(K_2\times P_n)$ and $D(K_2\times P_n)$. Our results coincide with those of Vince[14].
Section 5 proposes two questions on $K_m\times C_n$, where $C_n$ is the cycle of $n$ vertices.

\vspace{6mm}
\section{Notations and Preliminaries.}

By the definition of $K_m\times P_n$, 
it consists of $n$ subgraphs, each of which is isomorphic to $K_m$. For convenience, we refer to each such subgraph as a \emph{layer}. As shown in Figure 1, $H = K_3 \times P_3$ has three layers, where the cycle $v_{1,1}v_{1,2}v_{1,3}v_{1,1}$ forms the initial layer and $v_{3,1}v_{3,2}v_{3,3}v_{3,1}$ forms the third layer.

\begin{figure}[htbp]
\centering
\includegraphics[height=4cm]{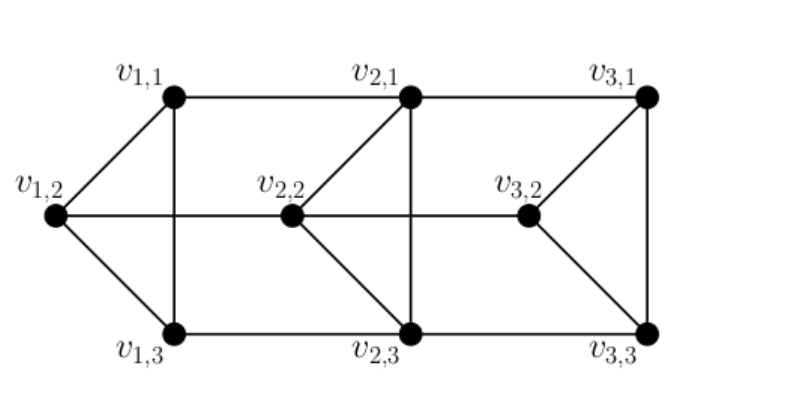}
\caption{$H$.}
\end{figure}

Let $f(m,k)$(where $1\leq k\leq n$) denote the number of connected sets in $K_m\times P_n$ which contain at least one vertex of each layer from the first to the $k$-th layer of $K_m\times P_n$. 
The following theorem follows directly.

\vspace{2mm}
{\bf Theorem 2.1[10].}\quad \emph{$N(K_m \times P_n)= \sum\limits_{k=1}^{n}(n+1-k)f(m,k)$.}\par
\vspace{4mm}

To compute $f(m,k)$, we partition it into $2^m-1$ parts. This is because every layer of $G$ contains $m$ vertices, and the connected sets that contain at least one vertex from each of the first 
$k$ layers induce a non-empty vertex subset in the $k$-th layer. There are $\binom{m}{1}+\binom{m}{2}+\dots+\binom{m}{m}(=2^{m}-1)$ such nonempty vertex subsets in the $k$-th layer. Denote these subsets by $\{S^k_{p}\}$ for $1\leq p \leq 2^{m}-1$.

Let $f_{S_{p}}(m,k)$ denote the number of connected subsets that contain at least one vertex of each layer from the first to the $(k-1)$-th layer
with exactly the vertices of $S^k_{p}$ in the $k$-th layer. Then 
$$f(m,k) = \sum_{p=1}^{2^m-1}f_{S_p}(m,k).$$

\vspace{2mm}
{\bf Lemma 2.2.}\quad \emph{For $1\leq p,q \leq 2^{m}-1$ with $p\neq q$ , if $|S^k_{p}|=|S^k_{q}|$, then $f_{S_{p}}(m,k)=f_{S_{q}}(m,k)$.}\par
\vspace{2mm}

The statement follows directly from symmetry.

\vspace{2mm}
We now introduce the following definitions which will be used in the formula for $A(K_m \times P_n)$. 
Let $f_i(m,k)$ be the number of connected sets in $K_m\times P_n$ which contains at least one vertex of each of the first $k-1$ layer and $i$ fixed vertices in the $k$-th layer, where $1\leq i \leq m$.  By Lemma 2.2,  we have $$f(m,k)=\sum_{i=1}^{m}\binom{m}{i}f_i(m,k). \eqno (2.1)$$
\vspace{2mm}

The following recurrence is taken from [10].

For $i=1,2,\dots, m$,

$$f_i(m,k)=\sum\limits_{j=1}^{m}\binom{m}{j}f_j(m,k-1)-\sum\limits_{j=1}^{m-i}\binom{m-i}{j}f_j(m,k-1).$$

Since $\binom{m-i}{j}=0$ for $j\geq m-i+1$,the recurrence for each $f_i(m,k)$ can simplify to:
$$f_i(m,k)=\sum_{j=1}^{m}[\binom{m}{j}-\binom{m-i}{j}]f_j(m,k-1).$$

\vspace{2mm}

Let

$$A_m = \left[
\begin{matrix}
\binom{m}{1}-\binom{m-1}{1}& \binom{m}{2}-\binom{m-1}{2}& \dots & \binom{m}{m} \\
\vdots & \vdots & \ddots &\vdots\\
\binom{m}{1}-\binom{m-i}{1}& \binom{m}{2}-\binom{m-i}{2}& \dots & \binom{m}{m}\\
\vdots & \vdots & \ddots &\vdots\\
\binom{m}{1} & \binom{m}{2}& \dots & \binom{m}{m}\end{matrix}
\right].
$$

Then $$[f_1(m,k), f_2(m,k), \dots, f_m(m,k)]^{T}=A_m[f_1(m,k-1), f_2(m,k-1), \dots, f_m(m,k-1)]^{T}.$$

Thus

$$[f_1(m,k), \dots, f_m(m,k)]^{T}=A_m[f_1(m,k-1),\dots, f_m(m,k-1)]^{T}$$
$$=A^2_m[f_1(m,k-2),\dots, f_m(m,k-2)]^{T}=\dots =A^k_m[1]^T, \eqno (2.2)$$

where $[1]^{T}$ denote the column vector $[1,1, \dots ,1]^{T}$ of length $m$. 

In particular,
$$f_m(m,k)= \sum_{i=1}^{m}\binom{m}{i}f_i(m,k-1)= f(m,k-1)$$ 
by the equation (2.1). For this reason, we call $A_m$ the \emph{recurrence matrix}.

\vspace{3mm}
Let $\mathbf{v}(i,k)$ denote the $i$-th column of $A^{k-1}_m$. Let $(A_m)^T$ denote the transpose of $A_m$. The following theorem establishes the relationship between $\mathbf{v}(i,k)$ and $f_i(m,k)$.

\vspace{2mm}
{\bf Theorem 2.3.}\quad \emph{$$[\binom{m}{1},\binom{m}{2},\dots,\binom{m}{m}]\times \mathbf{v}(i,k) = \binom{m}{i}f_i(m,k).$$}\par

\vspace{2mm} 
{\bf Proof.}
Let 
$$C_m = diag(\binom{m}{1},\binom{m}{2},\dots,\binom{m}{m})=\left[
\begin{matrix}
\binom{m}{1}& 0& \dots & 0 \\
0& \binom{m}{2}& \ddots & \vdots \\
\vdots& \ddots & \ddots &0\\
0&  \dots& 0 &\binom{m}{m} \end{matrix}
\right]
.$$

We first state a useful claim.

\vspace{2mm}
{\bf Claim 1.} For every $k\geq 1$, 
$$(C_m\times A^k_m)^T=C_m\times A^k_m.$$

\vspace{2mm}
{\bf Proof of Claim 1.}
We proceed by induction on $k$. The base case is $k=1$.
Using the identity $$\binom{m}{i}[\binom{m}{j}-\binom{m-i}{j}]=\binom{m}{i}\binom{m}{j}-\frac{m!}{i!j!(m-j-i)!}=\binom{m}{j}[\binom{m}{i}-\binom{m-j}{i}],$$

we observe that the $(i,j)$-entry of $C_mA_m$ equals the $(i,j)$-entry of $(A_m)^TC_m$.

Explicitly,
\begin{align*}
C_m\times A_m & = \left[
\begin{matrix}
\binom{m}{1}[\binom{m}{1}-\binom{m-1}{1}]& \dots & \binom{m}{1}[\binom{m}{i}-\binom{m-1}{i}]& \dots & \binom{m}{1}\binom{m}{m} \\
\vdots & \ddots&\vdots &\ddots &\vdots\\
\binom{m}{i}[\binom{m}{1}-\binom{m-i}{1}]&  \dots& \binom{m}{i}[\binom{m}{i}-\binom{m-i}{i}]& \dots & \binom{m}{i}\binom{m}{m}\\
\vdots & \ddots&\vdots &\ddots &\vdots\\
\binom{m}{m}\binom{m}{1} & \dots& \binom{m}{m}\binom{m}{i}& \dots  & \binom{m}{m}\binom{m}{m}\end{matrix}
\right]=(A_m)^T\times C_m.
\end{align*}

Assume $(C_mA^{k-1}_m)^T=C_mA^{k-1}_m$. Then 
$$(C_m\times A^k_m)^T=(A_m)^T\times (C_m\times A^{k-1}_m)^T=(A_m)^T\times C_m\times A^{k-1}_m=C_m\times A^k_m,$$
where we used the identity $C_m\times A_m=(A_m)^T\times C_m$ from the base case. This completes the induction.
$\hfill\qedsymbol$\par

\vspace{2mm}
Now return to the proof of Theorem 2.3.

By Claim 1, the matrix $C_mA^{k-1}_m$ is symmetric. Hence its row sums equal its column sums. The vector of row sums is 
$$[1,1,\dots, 1]\times C_mA^{k-1}_m=[\binom{m}{1}, \binom{m}{2}, \dots ,\binom{m}{m}]A^{k-1}_m.$$
On the other hand, the vector of column sums is 
$$C_mA^{k-1}_m\times [1]^T=[\binom{m}{1}f_1(m,k), \binom{m}{2}f_2(m,k), \dots, \binom{m}{m}f_m(m,k)]^T.$$
Since $A^{k-1}_m=[\mathbf{v}(1,k), \mathbf{v}(2,k),\dots,\mathbf{v}(i,k), \dots, \mathbf{v}(m,k)]$, comparing the $i$-th entries of the two equal vectors gives
$[\binom{m}{1},\binom{m}{2},\dots,\binom{m}{m}]\times \mathbf{v}(i,k) = \binom{m}{i}f_i(m,k)$.
$\hfill\qedsymbol$\par
\vspace{2mm}

\vspace{2mm}

Need to say that we can obtain some by-products using the recurrence matrix.
Now, let $p_m(\lambda)=|\lambda E_m - A_m|$, which is the characteristic polynomial of $A_m$. Then  
$$p_m(\lambda)=\lambda^{m}+c_{m,m}\lambda^{m-1}+ \dots + c_{m,1}. \eqno (2.3)$$ 

Since the value of $c_{m,m}$ is related to the Fibonacci sequence, let us recall some concepts about the Fibonacci sequence.
We denote the $n$-th term of the Fibonacci sequence by $F(n)$, where $F(n)=\frac{1}{\sqrt5}[(\frac{1+\sqrt5}{2})^{n}-(\frac{1-\sqrt5}{2})^{n}]$ and $F(0)=0, F(1)=1, F(2)=1, F(3)=2$.

\vspace{2mm}
{\bf Theorem 2.4.}\quad \emph{If $m\geq 3$, then $c_{m,m}= F(m+1)-2^{m}; c_{m,1} =1$. In particular, $c_{2,2}=-2, c_{2,1}=-1$.}\par

\vspace{2mm} 
{\bf Proof.} It is a well-known result that $-c_{m,m}$ equals the trail of the recurrence matrix $A_m$. Hence,
$$-c_{m,m}=\sum\limits^{m}_{i=1}\binom{m}{i}-\sum\limits^{\lfloor\frac{m}{2}\rfloor}_{i=1}\binom{m-i}{i}.$$
Obviously, the first part equals to $2^{m}-1$.
By [7], the second part is $\sum\limits^{\lfloor\frac{m}{2}\rfloor}_{i=0}\binom{m-i}{i} - \binom{m}{0}=F(m+1) - 1$. Therefore, we have $c_{m,m}= F(m+1)-1-2^{m} + 1= F(m+1)-2^m$.

The coefficient $c_{m,1}$ equals $(-1)^m|A_m|$. Every element of the last column of $A_m$ is $\binom{m}{m}=1$. 
Multiplying this column by $\binom{m}{i}$ and subtracting the result from the $i$-th column for each $1\leq i \leq m-1$, we obtain the following determinant.
$$
\begin{vmatrix}
-\binom{m-1}{1}& -\binom{m-1}{2}& \dots & -1 & 1 \\
-\binom{m-2}{1} & \vdots & -1 &0 & 1\\
\vdots & \vdots & \ddots &\vdots\\
-2& -1& \dots &  0& 1\\
-1& 0& \dots  &  0 &1\\
0 & 0& \dots  & 0  & 1
\end{vmatrix}
,
$$

in which the upper-left $(m-1)\times (m-1)$ block can be reduced to an upper-triangular matrix with -1 on every diagonal entry.
Since the reduction process involves an odd number of column exchanges (except when $m=2$), we obtain: $$ c_{m,1}= (-1)^m|A_m|= (-1)^m\times (-1)\times (-1)^{m-1} =(-1)^{2m}=1.$$ 
$\hfill\qedsymbol$\par
\vspace{2mm}

Consequently, we obtain the following linear recurrence:

\vspace{2mm}
{\bf Theorem 2.5.}\quad \emph{For every $k\geq m$, $f(m,k) = - c_{m,m} f(m,k-1) - c_{m,m-1} f(m,k-2) - \dots - c_{m,1}f(m,k-m)$.}\par

\vspace{2mm} 
{\bf Proof.} 
Let $\mathbf{0}_{m,m}$ denote the $m\times m$ zero matrix.
By the Cayley-Hamilton Theorem[13], we have $p_m(A_m)=\mathbf{0}_{m,m}$. It implies that $A_m^{m}+c_{m,m}A_m^{m-1}+ \dots + c_{m,1}E_m=\mathbf{0}_{m,m}$ by (2.3), where $E_m$ denotes the identity matrix of order $m$. 
Hence $$A_m^{m}=-c_{m,m}A_m^{m-1}- \dots - c_{m,1}E_m.$$

Multiplying both sides by $A^{k-m-1}_m\times [1]^{T}$ gives 

\begin{align*}
A_m^{k-1}\times [1]^T&=-c_{m,m}A_m^{k-2}\times [1]^T- \dots - c_{m,1}A^{k-m-1}_m\times [1]^T\\
[f_1(m,k),\dots,f_m(m,k)]^T&=-c_{m,m}[f_1(m,k-1),\dots,f_m(m,k-1)]^T- \dots\\ 
&- c_{m,1}[f_1(m,k-m),\dots,f_m(m,k-m)]^T.
\end{align*}

Multiplying both sides of the above equation by $[\binom{m}{1},\binom{m}{2},\dots,\binom{m}{m}]$ yields
\begin{align*}
\sum\limits_{i=1}^{m}\binom{m}{i}f_i(m,k)&=-c_{m,m}\sum\limits_{i=1}^{m}\binom{m}{i}f_i(m,k-1)-\dots-c_{m,1}\sum\limits_{i=1}^{m}\binom{m}{i}f_i(m,k-m-1).
\end{align*}

Applying (2.1) to each term, we obtain
\begin{align*}
f(m,k) &= - c_{m,m} f(m,k-1) - c_{m,m-1} f(m,k-2) - \dots - c_{m,1}f(m,k-m-1).
\end{align*}
$\hfill\qedsymbol$\par
\vspace{2mm}

\vspace{2mm}
{\bf Remark 2.6.}\quad \emph{Theorem 2.3 will be used in the proof of Lemma 3.4 in Section 3. Theorems 2.4 and 2.5 will be used in Section 4.}\par

\vspace{6mm}
\section{The average order of a connected set of $K_m\times P_n$.}

In this section, we discuss the average order of connected sets of $K_m\times P_n$. For $k=1,2,\dots, n$, let $F(m,k)$ denote the collection of connected sets of $K_m\times P_k$ that contain at least one vertex of each layer, and let $S(F(m,k))$ denote the sum of orders of all connected sets in $F(m,k)$. 

\vspace{2mm}
{\bf Theorem 3.1.} \emph{$$A(K_m\times P_n) = \frac{\sum\limits_{k=1}^{n}(n-k+1)S(F(m,k))}{\sum\limits_{k=1}^{n}(n-k+1)f(m,k)}.$$}
\vspace{2mm}

{\bf Proof.}
The numerator is precisely the total sum of the orders of all non-empty connected sets of $K_m\times P_n$, while the denominator equals $N(K_m\times P_n)$ by Theorem 2.1, the total number of such sets. Hence the formula gives the average order. $\hfill\qedsymbol$\par
\vspace{4mm}

Next, we intend to give a formula for $S(F(m,k))$.

For $p=1,2,\dots, 2^m-1$ and $k=1,2,\dots, n$, let $F_{S_p}(m,k)$ denote the collection of connected subsets in $K_{m}\times P_{n}$
that contain at least one vertex from each layer from the first to the $k-1$-th layer and exactly all vertices of $S^k_p$ in the $k$-th layer.

\vspace{2mm}
{\bf Lemma 3.2.}\quad \emph{For $1\leq p,q \leq 2^{m}-1$ with $p\neq q$ , if $|S^k_{p}|=|S^k_{p}|$, then the total sum of the orders of connected subsets in $F_{S_{p}}(m,k)$ equals that in $F_{S_{q}}(m,k)$.}\par
\vspace{2mm}

The statement follows directly from symmetry.

Now, for $i=1,2,\dots, m$ and $k=1,2,\dots, n$, let $F_i(m,k)$ denote the collection of connected sets of $K_m\times P_k$ that contain at least one vertex of each layer with exactly $i$ fixed vertices in the $k$-th layer. Let $s_{i}(m,k)$ denote the total sum of the orders
of connected sets in $F_i(m,k)$. 
Since every connected set of $K_m\times P_1$ has order $i$, we have $$[s_{1}(m,1), s_{2}(m,1), \dots ,s_{m}(m,1)]^{T}=[1, 2, \dots, m]^T.$$

By Lemma 3.2, $$S(F(m,k))=\sum_{i=1}^{m}\binom{m}{i}s_i(m,k).$$

To compute $s_i(m,k)$, we introduce a new matrix $B_m$ and use a recursive multiplication scheme involving matrices $B_m$ and $A_m$.
Let $[1]^{T}$ denote the column vector $[1,1, \dots ,1]^{T}$ of length $m$. Let 
$$B_m = diag(1,2,\dots,m)=\left[
\begin{matrix}
1& 0& \dots & 0 \\
0& 2& \ddots & \vdots \\
\vdots& \ddots & \ddots &0\\
0&  \dots& 0 &m \end{matrix}
\right]
.$$

Let $E_m$ denote the identity matrix of order $m$. Then we have the following lemma.

\vspace{2mm}
{\bf Lemma 3.3.} \emph{$$[s_{1}(m,k), s_{2}(m,k), \dots ,s_{m}(m,k)]^{T}=(\sum\limits^{k-1}_{s=0}A^{k-s-1}_mB_mA^{s}_m) \times [1]^{T},$$ where $A^0_m=E_m$.}
\vspace{4mm}

{\bf Proof.} We use induction on $k$. The base case is that $k=1$. Since $f_{i}(m,1) = 1$ for $i=1,2,\dots, m$, we have
\begin{align*}
  [s_{1}(m,1), s_{2}(m,1), \dots, s_{m}(m,1)]^{T}&=B_m\times [1]^{T}\\
  & =[1f_{1}(m,1),2f_{2}(m,1), \dots, mf_{m}(m,1)]^{T}.
\end{align*}

So the statement holds for $k=1$. Assume the statement true for $k-1$ and consider the case $k$. Let $$[s'_{1}(m,k), s'_{2}(m,k), \dots ,s'_{m}(m,k)]^{T}=A_m[s_{1}(m,k-1), s_{2}(m,k-1), \dots ,s_{m}(m,k-1)]^{T}.$$ By the definition of $A_m$ and the induction hypothesis, $s'_{i}(m,k)$ equals the sum of numbers of vertices in $F_i(m,k)$ that lie in layers other than the $k$-th layer.

On the other hand, the vector $$[1f_{1}(m,k),2f_{2}(m,k), \dots, mf_{m}(m,k)]^{T}$$ gives the sum of numbers of vertices in $F_{i}(m,k)$ that lie in the $k$-th layer for $i=1,2,\dots, m$. By the equation (2.2), $$[f_{1}(m,k), f_{2}(m,k), \dots ,f_{m}(m,k)]^{T}= A^{k-1}_m\times [1]^{T}.$$ Hence $$[1f_{1}(m,k),2f_{2}(m,k), \dots, mf_{m}(m,k)]^{T}=B_m\times [f_{1}(m,k), f_{2}(m,k), \dots ,f_{m}(m,k)]^{T}=B_mA^{k-1}_m\times [1]^{T}.$$

Therefore, $$[s_{1}(m,k), s_{2}(m,k), \dots ,s_{m}(m,k)]^{T}$$
\begin{align*}
&=[s'_{1}(m,k), s'_{2}(m,k), \dots ,s'_{m}(m,k)]^{T}+[1f_{1}(m,k),2f_{2}(m,k), \dots, mf_{m}(m,k)]^{T}\\
   &= A_m\times[s_{1}(m,k-1), s_{2}(m,k-1), \dots ,s_{m}(m,k-1)]^{T}\times [1]^{T} + B_m\times A^{k-1}_m\times [1]^{T} \\
   &= A_m\times \sum\limits^{k-2}_{s=0}A^{k-s-2}_mB_mA^{s}_m \times [1]^{T} + E_m\times B_m\times A^{k-1}_m\times [1]^{T} \\
   &= \sum\limits^{k-1}_{s=0}A_m^{k-s-1}B_mA^{s}_m\times [1]^{T},
\end{align*}
which completes the induction.$\hfill\qedsymbol$\par
\vspace{2mm}

However, obtaining $s_i(m,k)$ requires evaluating the matrix expression in Lemma 3.3. To simplify the computation we introduce an auxiliary result.

Let $E_{i,i}$ be the $m\times m$ matrix whose $(i,i)$-entry is 1 and all other entries are 0. We note that $B_m = mE_m-(m-1)E_{1,1}-(m-2)E_{2,2}-\dots -E_{m-1,m-1}$.
In order to compute $\sum\limits_{s=0}^{k-1}A^{k-s-1}_mB_mA^s_m \times [1]^T$, we first calculate $\sum\limits_{s=0}^{k-1}A^{k-s-1}_mE_{i,i}A^s_m \times [1]^T$ for $i=1,2,\dots, m$.

For each $j=1,2,\dots,m$, we define 
$$[x_{i,1}(k), x_{i,2}(k), \dots, x_{i,j}(k), \dots, x_{i,m}(k)]^{T} = \sum\limits^{k-1}_{s=0}A^{k-s-1}_mE_{i,i}A^{s}_m\times [1]^{T}.$$
The values of $s_i(m,k)$ can then be obtained from $x_{i,j}(k)$.

\vspace{2mm}
{\bf Lemma 3.4.} \emph{$\sum\limits^{m}_{j=1}\binom{m}{j} x_{i,j}(k)=\sum\limits^{k}_{s=1}\binom{m}{i}f_{i}(m,s)f_{i}(m,k+1-s).$}
\vspace{2mm}

{\bf Proof.} Since 

$$E_{i,i} =  \left[
\begin{matrix}
0 & \cdots & 0 & \cdots & 0 \\
\vdots & \ddots & \vdots & & \vdots \\
0 & \cdots & 1 & \cdots & 0 \\
\vdots & & \vdots & \ddots & \vdots \\
0 & \cdots & 0 & \cdots & 0
\end{matrix} \right],
$$

we have 

$$
A_m^{k-s-1} E_{i,i} A_m^{s}
= \left[\begin{matrix}
a_{1,i}b_{i,1} & a_{1,i}b_{i,2} & \cdots & a_{1,i}b_{i,m} \\
a_{2,i}b_{i,1} & a_{2,i}b_{i,2} & \cdots & a_{2,i}b_{i,m} \\
\vdots & \vdots & \ddots & \vdots \\
a_{m,i}b_{i,1} & a_{m,i}b_{i,2} & \cdots & a_{m,i}b_{i,m}
\end{matrix}\right],
$$
where $a_{j,i}$ is the $(j,i)$-entry of $A^{k-s-1}_m$ and $b_{i,j}$ is the $(i,j)$-entry of $A^{s}_m$.

Because $A^{s}_m\times [1]^{T} = [f_1(m,s+1), f_2(m,s+1), \dots , f_m(m,s+1)]^{T}$, we have $\sum\limits^{m}_{j=1}b_{i,j}=f_{i}(m,s+1)$.
Hence the sum of the $j$-th row of $A^{k-s-1}_mE_{i,i}A^{s}_m$ is $a_{j,i}f_{i}(m,s+1)$.

Since $[a_{1,i}, a_{2,i}, \dots, a_{m,i}]^T$ is the $i$-th column of $A^{k-s-1}_m$, Theorem 2.3 implies $$[\binom{m}{1},\binom{m}{2},\dots, \binom{m}{m}]\times [a_{1,i}, a_{2,i}, \dots, a_{m,i}]^T=\binom{m}{i}f_{i}(m,k-s).$$

Consequently,
\begin{align*}
[x_{i,1}(k),x_{i,2}(k),\dots,x_{i,m}(k)]^T=\sum\limits^{k-1}_{s=0}[a_{1,i}f_i(m,s+1),a_{2,i}f_i(m,s+1),\dots, a_{m,i}f_i(m,s+1)]^{T}
\end{align*}

Multiplying both sides of the above equation by $[\binom{m}{1},\binom{m}{2},\dots,\binom{m}{m}]$ yields
\begin{align*}
\sum\limits^{m}_{j=1}\binom{m}{j} x_{i,j}(k) &=[\binom{m}{1},\dots,\binom{m}{m}]\times \sum\limits^{k-1}_{s=0}[a_{1,i}f_i(m,s+1),\dots, a_{m,i}f_i(m,s+1)]^{T} \\
&=\sum\limits^{k-1}_{s=0}\binom{m}{i}f_{i}(m,k-s)f_{i}(m,s+1)\\
& =\sum\limits_{s=1}^{k}\binom{m}{i}f_{i}(m,s)f_{i}(m,k+1-s).
\end{align*}

The lemma holds.
$\hfill\qedsymbol$\par
\vspace{2mm}

\vspace{2mm}
{\bf Theorem 3.5.} \emph{$S(F(m,k))=mkf(m,k)-\sum\limits_{i=1}^{m-1}\sum\limits_{s=1}^{k}\binom{m}{i}(m-i)f_{i}(m,s)f_{i}(m,k+1-s).$}
\vspace{2mm}

{\bf Proof.}
Observe that 
\begin{align*}
B_m & =mE_m- \sum\limits^{m-1}_{i=1}(m-i)E_{i,i}.
\end{align*}

Then \begin{align*}
S(F(m,k))&= \sum\limits^{m}_{i=1}\binom{m}{i} s_i(m,k)\\
&=[\binom{m}{1},\binom{m}{2},\dots,\binom{m}{m}]\times[s_1(m,k),s_2(m,k),\dots,s_m(m,k)]^T\\
&=[\binom{m}{1},\binom{m}{2},\dots,\binom{m}{m}]\times\sum\limits^{k-1}_{s=0}A^{k-s-1}B_mA^{s}\times [1]^{T}.\\
\end{align*}

Applying Lemmas 3.4, we obtain
\begin{align*}
S(F(m,k)) &=[\binom{m}{1},\binom{m}{2},\dots,\binom{m}{m}]\times [mkE_mA^{k-1}- \sum\limits_{i=1}^{m-1}(m-i)\sum\limits^{k-1}_{s=0}A^{k-s-1}E_{i,i}A^{s}]\times [1]^{T}\\
&=mk[\binom{m}{1},\binom{m}{2},\dots,\binom{m}{m}]\times[f_1(m,k),f_2(m,k),\dots,f_m(m,k)]^T\\
&-  \sum\limits_{i=1}^{m-1}(m-i)[\binom{m}{1},\binom{m}{2},\dots,\binom{m}{m}]\times [x_{i,1}(k), x_{i,2}(k), \dots, x_{i,j}(k), \dots, x_{i,m}(k)]^{T}\\
&=\sum\limits^{m}_{i=1}\binom{m}{i}f_i(m,k) -\sum\limits_{i=1}^{m-1}\sum\limits^{m}_{j=1}\binom{m}{j}(m-i)x_{i,j}(k)\\
&=mkf(m,k)-\sum\limits_{i=1}^{m-1}\sum\limits_{s=1}^{k}\binom{m}{i} (m-i)f_{i}(m,s)f_{i}(m,k+1-s),
\end{align*}

which completes the proof.
$\hfill\qedsymbol$\par
\vspace{2mm}

Combining Theorems 3.1 and 3.5, we have the general formula for the average order and the density of connected sets in $K_m\times P_n$. 

\vspace{2mm}
{\bf Theorem 3.6.} \emph{$$A(K_m\times P_n) = \frac{\sum\limits_{k=1}^{n}(n-k+1)\{mkf(m,k)-\sum\limits_{i=1}^{m-1}\sum\limits_{s=1}^{k}\binom{m}{i}(m-i)f_{i}(m,s)f_{i}(m,k+1-s)\}}{\sum\limits_{k=1}^{n}(n-k+1)f(m,k)}.$$}

\vspace{2mm}
{\bf Theorem 3.7.} \emph{$$D(K_m\times P_n) = \frac{\sum\limits_{k=1}^{n}(n-k+1)\{mkf(m,k)-\sum\limits_{i=1}^{m-1}\sum\limits_{s=1}^{k}\binom{m}{i}(m-i)f_{i}(m,s)f_{i}(m,k+1-s)\}}{mn\sum\limits_{k=1}^{n}(n-k+1)f(m,k)}.$$}
\vspace{4mm}

\section{The number and average order of connected sets in $K_2\times P_n$.}
In this section, we apply the matrix method developed in Sections 2 and 3 to study the number and the average order of connected sets in $K_2\times P_n$.
Need to say that Vince[14] has given a formula for $A(K_2\times P_n)$(cf. formula(1.1)). However, our method is different from that in [14].

For $m=2$, the recurrence matrix is 
 $$A_2 = \left[
\begin{matrix}
1& 1 \\
2& 1 \end{matrix}
\right].
$$
By Theorem 2.4, the characteristic polynomial of $A_2$ is $p_2(\lambda) = \lambda^2 -2 \lambda - 1$. Solving $p_2(\lambda) = 0$ gives the two eigenvalues: $$\lambda_1 = 1-\sqrt{2}, \quad \lambda_2=1+\sqrt{2}.$$

By Theorem 2.5, we obtain the recurrence $$f(2,k) = 2 f(2,k-1) + f(2,k-2), \qquad k \geq 2,$$ 
with initial values $f(2,1)=3$ and $f(2,2)=7$.

\vspace{2mm}
{\bf Lemma 4.1.}\quad \emph{$$f(2,k)=\frac{(1-\sqrt{2})^{k+1}}{2}+\frac{(1+\sqrt{2})^{k+1}}{2}.$$}\par
\vspace{2mm}

\vspace{2mm}
{\bf Proof.}
Assume that $f(2,k)=C_1(1-\sqrt{2})^k+C_2(1+\sqrt{2})^k$. 
Using the initial values $f(2,1)=3$ and $f(2,2)=7$, we obtain: $C_1=\frac{1-\sqrt{2}}{2}$ and $C_2=\frac{1+\sqrt{2}}{2}$.
Substituting these coefficients gives the stated formula.
$\hfill\qedsymbol$\par
\vspace{2mm}

One may note that $f(2,k)=\beta(k+1)$, where $\beta(k)$ denotes the $k$-th Pell-Lucas number[7]. 

\vspace{2mm}
{\bf Theorem 4.2.}\quad \emph{$$N(K_2 \times P_n)=\frac{f(2,n+2)-4n-7}{2}.$$}\par
\vspace{2mm}

\vspace{2mm}
{\bf Proof.} By Theorem 2.1, $$N(K_2 \times P_n)= \sum\limits_{k=1}^{n}(n+1-k)f(2,k) =(n+1)\sum\limits_{k=1}^{n}f(2,k) -\sum\limits_{k=1}^{n}kf(2,k).$$
Applying formulas given in [7], we have 

$$\sum\limits_{k=1}^{n}f(2,k)=\frac{f(2,n+1)+f(2,n)-4}{2}.\eqno(4.1)$$

Next we compute $\sum\limits_{k=1}^{n}kf(2,k)$.
By Lemma 4.1, $$\sum\limits^{n}_{k=1}kf(2,k)=\sum\limits^{n}_{k=1}\frac{k(1-\sqrt{2})^{k+1}}{2}+\sum\limits^{n}_{k=1}\frac{k(1+\sqrt{2})^{k+1}}{2} \eqno(4.2)$$

Since $\sum\limits^{n}_{k=1}kx^{k+1}=x^2(\sum\limits^{n}_{k=1}x^{k})^{'}$ and $(\sum\limits^{n}_{k=1}x^{k})^{'}=\frac{1+nx^{n+1}-(n+1)x^{n}}{(1-x)^{2}}$, the right part of the equation (4.2) can be transfer to
$$\frac{(1-\sqrt{2})^2+n(1-\sqrt{2})^{n+3}-(n+1)(1-\sqrt{2})^{n+2}}{4}+\frac{(1+\sqrt{2})^2+n(1+\sqrt{2})^{n+3}-(n+1)(1+\sqrt{2})^{n+2}}{4}.$$

After simplification this equals $$\frac{2nf(2,n+2)-(2n+2)f(2,n+1)+6}{4},$$

which gives 
$$\sum\limits_{k=1}^{n}kf(2,k) =\frac{nf(2,n+2)-(n+1)f(2,n+1)+3}{2}. \eqno(4.3)$$

Now substitute (4.1) and (4.3) into the expression for $N(K_2 \times P_n)$: 
\begin{align*}
  N(K_2 \times P_n) & = \frac{(n+1)f(2,n+1)+(n+1)f(2,n)-4(n+1)}{2}-\frac{nf(2,n+2)-(n+1)f(2,n+1)+3}{2}\\
  & =\frac{(2n+2)f(2,n+1)+(n+1)f(2,n)f(2,n)-nf(2,n+2)-4n-7}{2}\\
  &=\frac{f(2,n+2)-4n-7}{2},
\end{align*}
where the last equality follows from the recurrence $f(2,n+2)=2f(2,n+1)+f(2,n)$.
$\hfill\qedsymbol$\par
\vspace{2mm}

\vspace{2mm}
We now turn to $A(K_2\times P_n)$. By Theorem 3.1, this requires the expression of $S(F(2,k))$.

By Theorem 3.5, 
$$S(F(2,k)) = 2kf(2,k)-2\sum\limits^{k}_{i=1}f_{1}(2,i)f_{1}(2,k+1-i).$$

\vspace{2mm}
{\bf Lemma 4.3.}\quad \emph{$$f_1(2,k)=\frac{(1+\sqrt{2})^{k}}{2\sqrt{2}}-\frac{(1-\sqrt{2})^{k}}{2\sqrt{2}}.$$}\par
\vspace{2mm}

\vspace{2mm}
{\bf Proof.}
Assume that $f_1(2,k)=C_3(1-\sqrt{2})^k+C_4(1+\sqrt{2})^k$. Multiplying $A_2$ and $A^2_2$ by $[1,1]^T$ gives the initial values $f_1(2,1)=1$ and $f_1(2,2)=2$. Solving for the coefficients yields: $C_3=\frac{-1}{2\sqrt{2}}$ and $C_4=\frac{1}{2\sqrt{2}}$.
$\hfill\qedsymbol$\par
\vspace{2mm}

Lemma 4.3 shows that $f_{1}(2,k)=\bar{\beta}(k)$, where $\bar{\beta}(k)$ is the $k$-th Pell number[7].
From [7, p.212], we know that $\sum\limits^{k}_{i=1}f_{1}(2,i)f_{1}(2,k+1-i) = \frac{(k+2)f(2,k) - f_{1}(2,k+2)}{4}$.
Therefore
\begin{align*}
  S(F(2,k)) & =2kf(2,k)- \frac{(k+2)f(2,k) - f_{1}(2,k+2)}{4} \times 2  \\
  & = \frac{(3k-2)f(2,k) + f_1(2,k+2)}{2}.
\end{align*}

Thus $$\sum\limits_{k=1}^{n}(n-k+1)S(F(k)) = \sum\limits_{k=1}^{n}(n-k+1)\frac{(3k-2)f(2,k) + f_1(2,k+2)}{2}.$$

Using the identities given in [7], we have 
\begin{align*}
\sum\limits_{k=1}^{n}f_1(2,k+2) &= \frac{f(2,n+3) - 7}{2}.
\end{align*}
 
Following the same method as in Theorem 4.2, we also obtain
\begin{align*}
\sum\limits_{k=1}^{n}kf_1(2,k+2) &= \frac{1}{2}[2(n-1)f_1(2,n+2)+(3n-1)f_1(2,n+1)+nf_1(2,n)+5],\\  
\sum\limits_{k=1}^{n}k^2f(2,k) &= \frac{1}{2}[(2n^2+2n+1)f_1(2,n+2)+(1-2n)f_1(2,n+3)-7].
\end{align*}

From these we derive

\begin{align*}
\sum\limits_{k=1}^{n}(n-k+1)S(F(k)) & = \frac{(21n-44)f(2,n)+(-12n+26)f_1(2,n)+17n+44}{4}\\
   &+\frac{12f(2,n)-7f_1(2,n)-7n-12}{4}\\
   &=\frac{(21n-32)f(2,n)+(-12n+19)f_1(2,n)+10n+32}{4}.
\end{align*}

Combining this with Theorem 4.2 gives the following closed-form expressions.

\vspace{2mm}
{\bf Theorem 4.4.}\quad \emph{$$A(K_2 \times P_n)= \frac{(21n-32)f(2,n)+(-12n+19)f_1(2,n)+10n+32}{2f(2,n+2)-4n-7}.$$}\par
\vspace{2mm}

\vspace{2mm}
{\bf Theorem 4.5.}\quad \emph{$$D(K_2 \times P_n)= \frac{(21n-32)f(2,n)+(-12n+19)f_1(2,n)+10n+32}{2n[2f(2,n+2)-4n-7]}.$$}\par
\vspace{2mm}

\vspace{2mm}
{\bf Remark 4.6.}\quad \emph{Since $f(2,n)=f(2,n-1)+2f_1(2,n)$,
$$\frac{(21n-32)f(2,n)+(-12n+19)f_1(2,n)+10n+32}{2f(2,n+2)-4n-7}= $$
$$\frac{[32-45f_1(2,n)- 32f(2,n-1)]+ n[10+21f(2,n-1)+30f_1(2,n)]}{4},$$
which shows that our formulas coincide with those obtained by Vince in [14].}\par

\vspace{6mm}

\section{Further works.}
We have discussed the number and average order of connected sets in $K_m\times P_n$. A natural follow-up question is whether analogous formulas can be obtained for $K_m\times C_n$.
Specifically:

\vspace{2mm}
{\bf Question 1.} Find a formula for $N(K_m\times C_n)$.

\vspace{4mm}
{\bf Question 2.} Find a formula for $A(K_m\times C_n)$.
\vspace{2mm}

Exploring these questions would extend the results presented here and provide a more complete picture of connected sets in products of complete graphs with simple linear structures.

\end{CJK*}

\end{spacing}
\end{document}